# Moment estimation for ergodic diffusion processes

YURY A. KUTOYANTS[1] and NAKAHIRO YOSHIDA[2]

[1]*Laboratoire de Statistique et Processus, Université du Maine, 72017 Le Mans, France.
E-mail: kutoyants@univ-lemans.fr*
[2]*Graduate School of Mathematical Sciences, The University of Tokyo, Komaba 3-8-1, Meguro-ku, Tokyo 153-0041, Japan. E-mail: nakahiro@ms.u-tokyo.ac.jp*

We investigate the moment estimation for an ergodic diffusion process with unknown trend coefficient. We consider nonparametric and parametric estimation. In each case, we present a lower bound for the risk and then construct an asymptotically efficient estimator of the moment type functional or of a parameter which has a one-to-one correspondence to such a functional. Next, we clarify a higher order property of the moment type estimator by the Edgeworth expansion of the distribution function.

## 1. Introduction

Suppose that a diffusion process $\{X_t, 0 \leq t \leq T\}$ is uniquely defined by the stochastic differential equation

$$\mathrm{d}X_t = S(X_t)\,\mathrm{d}t + \sigma(X_t)\,\mathrm{d}W_t, \qquad 0 \leq t \leq T, \tag{1}$$

where $\{W_t, t \geq 0\}$ is a standard Wiener process, and $X_0$ is an initial random variable independent of $\{W_t, t \geq 0\}$. We assume that $\{X_t, 0 \leq t \leq T\}$ is ergodic with the invariant probability distribution $\mu_S$ (depending on $S$ and $\sigma$).

In our statistical problem, the diffusion coefficient $\sigma$ is known to the observer, while the trend coefficient $S$ is unknown. Given a function $F: \mathbf{R} \to \mathbf{R}$, we want to estimate the parameter

$$\vartheta = \vartheta_S = \mathbf{E}_S[F(\xi)] = \int F(y)\mu_S(\mathrm{d}y),$$

where $\mathbf{E}_S$ denotes the expectation with respect to $\mu_S$, and $\xi$ denotes a "stationary" random variable with distribution $\mu_S$. If $X_t$ is stationary, we take $\xi = X_0$; if not, we may enlarge the original probability space to realize $\xi$.

We are interested in the asymptotically efficient estimation of $\vartheta$ based on the observations $\{X_t, 0 \leq t \leq T\}$ as $T \to \infty$. We will consider nonparametric and parametric estimation problems.

Section 2 treats a nonparametric case where the function $S$ is completely unknown. As usual in such problems, we derive a minimax bound of the risks of all estimators and then show that the *empirical moment estimator*

$$\vartheta_T^* = \frac{1}{T} \int_0^T F(X_t) \, dt$$

is asymptotically efficient in that $\vartheta_T^*$ attains this lower bound.

In Section 3, we suppose that the function $S$ belongs to a parametric family, that is, $S = S(\gamma, \cdot), \gamma \in \Gamma$, and $\vartheta$ is a function of $\gamma$: $\vartheta = \vartheta(\gamma)$. Then it turns out that the maximum likelihood estimator (MLE) $\hat{\gamma}_T$ (under certain regularity conditions) provides the asymptotically efficient estimator $\hat{\vartheta}_T = \vartheta(\hat{\gamma}_T)$ of $\vartheta$. However, the computation of the MLE in nonlinear models is often not easy. To avoid this drawback, we consider the so-called *one-step MLE* introduced by LeCam [6], which allows us to improve the empirical moment estimator to an asymptotically efficient one in this parametric setting. Here the asymptotic efficiency is in the sense of Hajek and LeCam.

After investigations of the first-order asymptotics, a natural direction of study is the second-order approximation of the distribution of the estimator. In Section 4, we derive an asymptotic expansion formula with the help of the *local approach* of the Malliavin calculus developed recently by Yoshida [12, 13, 14] and Kusuoka and Yoshida [2].

## 2. Nonparametric estimation

We consider the diffusion process $\{X_t, \, 0 \leq t \leq T\}$ described in Section 1.

In this section, $S$ is unknown, while $\sigma$ is a known continuous positive function. Let us denote by **S** the class of functions $S$ that satisfies

$$V(S, x) = \int_0^x \exp\left\{-2 \int_0^y \frac{S(v)}{\sigma(v)^2} \, dv\right\} dy \to \pm\infty \qquad \text{as } x \to \pm\infty \tag{2}$$

and

$$G(S) = \int_{-\infty}^{\infty} \sigma(y)^{-2} \exp\left\{2 \int_0^y \frac{S(v)}{\sigma(v)^2} \, dv\right\} dy < \infty. \tag{3}$$

These conditions guarantee the existence of invariant probability measure $\mu_S$ with the density function

$$f_S(x) = G(S)^{-1} \sigma(x)^{-2} \exp\left\{2 \int_0^x \frac{S(v)}{\sigma(v)^2} \, dv\right\} \tag{4}$$

and the law of large numbers. We suppose that the initial value $X_0 = \xi$ has a probability density function $f(\cdot)$, in particular, the process $X_t, \, t \geq 0$, is stationary.

We would like to estimate the mathematical expectation

$$\vartheta_S = \int_{-\infty}^{\infty} F(y) f_S(y)\, \mathrm{d}y.$$

To construct a lower minimax bound on the risks of all estimators, we define a nonparametric vicinity of a fixed model as follows. Fix some $S_*(\cdot) \in \mathbf{S}$ and $\delta > 0$, and introduce the set

$$V_\delta = \left\{ S(\cdot) : \sup_{x \in R} |S(x) - S_*(x)| < \delta \right\}.$$

We suppose that $S_*(\cdot)$ is such that the conditions (2) and (3) are fulfilled for all $S(\cdot) \in V_\delta$ and

$$\sup_{S(\cdot) \in V_\delta} G(S) < \infty, \qquad \sup_{S(\cdot) \in V_\delta} \mathbf{E}_S |F(\xi)| < \infty. \tag{5}$$

The role of Fisher information in our problem will be played by the quantity

$$\mathrm{I}(S) = \left\{ 4 \mathbf{E}_S \left[ \left( \frac{M_S(\xi)}{\sigma(\xi) f_S(\xi)} \right)^2 \right] \right\}^{-1}, \tag{6}$$

where $M_S(y) = \mathbf{E}_S[(F(\xi) - \vartheta_S)\chi_{\{\xi < y\}}]$. We put $\mathrm{I}_* = \mathrm{I}(S_*)$. We choose the polynomial loss function $\ell(u) = |u|^p$, $p > 0$.

**Theorem 1.** *Let the conditions (5) be fulfilled and let $\mathrm{I}_* > 0$. Then*

$$\lim_{\delta \to 0} \lim_{T \to \infty} \inf_{\bar\vartheta_T} \sup_{S(\cdot) \in V_\delta} \mathbf{E}_S[\ell(T^{1/2}(\bar\vartheta_T - \vartheta_S))] \geq \mathbf{E}[\ell(\eta \mathrm{I}_*^{-1/2})], \tag{7}$$

*where $\mathcal{L}(\eta) = \mathcal{N}(0,1)$ and inf is taken over all possible estimators $\bar\vartheta_T$ of the unknown parameter.*

**Proof.** We follow the well known scheme described in Ibragimov and Khasminskii ([1], Chapter 4); see also the proofs of similar results in Kutoyants [4].

Denote $\vartheta_* = \vartheta_{S_*}$ and let $\psi(\cdot)$ be a continuous function with compact support such that $S_h(\cdot) = S_*(\cdot) + (h - \vartheta_*)\psi(\cdot)\sigma(\cdot)^2 \in V_\delta$ for all $h \in (\vartheta_* - \gamma, \vartheta_* + \gamma)$. Here $\gamma > 0$ is a number chosen in such a way that $S_h(\cdot) \in V_\delta$. For the process (1) with $S(\cdot) = S_h(\cdot)$ we consider the parameter estimation problem for $h$ and recall the construction of the Hajek–LeCam minimax bound in this situation.

The direct expansion of the function $\vartheta_h = \vartheta_{S_* + (h - \vartheta)\psi(\cdot)\sigma(\cdot)^2}$ by the powers of $h - \vartheta_*$ gives the representation

$$\vartheta_h = \vartheta_* + 2(h - \vartheta_*)\left[ \mathbf{E}\left[ F(\xi) \int_0^\xi \psi(v)\, \mathrm{d}v \right] - \vartheta_* \mathbf{E}\left[ \int_0^\xi \psi(v)\, \mathrm{d}v \right] \right] + o(h - \vartheta_*),$$

where we write $\mathbf{E} = \mathbf{E}_{S_*}$ and $\mathbf{E}_h = \mathbf{E}_{S_h}$. Therefore, if we take $\psi(\cdot)$ from the class

$$\mathcal{K} = \left\{ \psi(\cdot) : \mathbf{E}\left[ [F(\xi) - \vartheta_*] \int_0^\xi \psi(v)\,\mathrm{d}v \right] = 2^{-1} \right\},$$

then $\vartheta_h = h + o(h - \vartheta_*)$. The corresponding family of measures $\{\mathbf{P}_h^{(T)}, h \in (\vartheta_* - \gamma, \vartheta_* + \gamma)\}$ is locally asymptotically normal (LAN) at the point $h = \vartheta_*$, therefore, we can apply the Hajek–LeCam inequality (see Ibragimov and Khasminski [1], Theorem 2.12.1) and write for $\psi \in \mathcal{K}$,

$$\lim_{\delta \to 0} \varliminf_{T \to \infty} \inf_{\bar{\vartheta}_T} \sup_{S(\cdot) \in V_\delta} \mathbf{E}_S[\ell(T^{1/2}(\bar{\vartheta}_T - \vartheta_S))]$$

$$\geq \lim_{\gamma \to 0} \varliminf_{T \to \infty} \inf_{\bar{\vartheta}_T} \sup_{|h - \vartheta| < \gamma} \mathbf{E}_h[\ell(T^{1/2}(\bar{\vartheta}_T - \vartheta_h))] \geq \mathbf{E}[\ell(\eta \mathrm{I}_\psi^{-1/2})].$$

Here $\eta \sim N(0,1)$ and $\mathrm{I}_\psi = \mathbf{E}(\psi(\xi)\sigma(\xi))^2$ is the Fisher information in the problem of estimation of $h$.

Furthermore, using the integration-by-parts formula and the Cauchy–Schwarz inequality we can write

$$2^{-1} = -\int_{-\infty}^{\infty} \int_{-\infty}^{y} (F(z) - \vartheta_*) f_{S_*}(z) \mathrm{d}z\, \psi(y)\, \mathrm{d}y$$

$$+ \left[ \int_{-\infty}^{y} (F(z) - \vartheta_*) f_{S_*}(z)\, \mathrm{d}z \int_0^y \psi(v)\, \mathrm{d}v \right]_{-\infty}^{\infty}$$

$$\leq \mathrm{I}_\psi^{1/2} \left\{ \mathbf{E}\left[ \left( \frac{M_{S_*}(\xi)}{\sigma(\xi) f_{S_*}(\xi)} \right)^2 \right] \right\}^{1/2}.$$

Therefore,

$$\mathrm{I}_\psi \geq \left\{ 4\mathbf{E}\left[ \left( \frac{M_{S_*}(\xi)}{\sigma(\xi) f_{S_*}(\xi)} \right)^2 \right] \right\}^{-1} \equiv \mathrm{I}_*.$$

We will have the equality if we choose

$$\psi(v) = \psi_*(v) \equiv \frac{C}{\sigma(v)^2 f_{S_*}(v)} \int_{-\infty}^{v} (F(y) - \vartheta_*) f_{S_*}(y)\, \mathrm{d}y$$

with the corresponding normalizing constant $C > 0$, but this function does not have compact support and, therefore, it cannot belong to $\mathcal{K}$. As usual in such situation (see Ibragimov and Khasminski [1], page 218), we introduce a sequence of smooth functions $\{\psi_N(\cdot)\}$ with compact supports approximating $\psi_*(\cdot)$ and such that

$$\inf_{\psi(\cdot) \in \mathcal{K}} \mathrm{I}_\psi = \lim_{N \to \infty} \mathrm{I}_{\psi_N} = \mathrm{I}_*.$$

Then

$$\sup_{\psi(\cdot)\in\mathcal{K}} \mathbf{E}\ell(\eta \mathrm{I}_\psi^{-1/2}) = \mathbf{E}\ell(\eta \mathrm{I}_*^{-1/2})$$

and we have the desired estimate (7). $\square$

**Definition 1.** *We say that the estimator $\bar{\vartheta}_T$ is asymptotically efficient for the loss function $\ell(\cdot)$ if*

$$\lim_{\delta\to 0}\lim_{T\to\infty}\sup_{S(\cdot)\in V_\delta} \mathbf{E}_S\ell(T^{1/2}(\bar{\vartheta}_T - \vartheta_S)) = \mathbf{E}\ell(\eta \mathrm{I}_*^{-1/2}) \qquad (8)$$

*for all functions $S_*(\cdot) \in \mathbf{S}$.*

Below we will show that the empirical estimator is an asymptotically efficient estimator of $\vartheta$ in this sense.

The process $\{X_t\}_{t\in\mathbf{R}_+}$ is stationary; therefore the estimator $\vartheta_T^*$ is unbiased: $\mathbf{E}_S\vartheta_T^* = \vartheta_S$. If the initial value $X_0$ has another distribution, then the estimator $\vartheta_T^*$ is no longer unbiased but nevertheless we have $|\mathbf{E}_S\vartheta_T^* - \vartheta_S| \leq \frac{C}{T}$ (see Kutoyants [3], Lemma 3.4.8) and the properties of the estimator $\vartheta_T^*$ can be studied without assumption of stationarity as well.

Let us introduce the function

$$H_S(y) = \int_0^y \frac{2}{\sigma(x)^2 f_S(x)} \int_{-\infty}^x [F(v) - \vartheta_S] f_S(v)\,dv\,dx$$

and the conditions that, for some $p_* > 0$,

$$\sup_{S(\cdot)\in V_\delta} \mathbf{E}_S|H_S(\xi)|^{p_*} < \infty, \qquad \sup_{S(\cdot)\in V_\delta} \mathbf{E}_S\left|\frac{M_S(\xi)}{\sigma(\xi)f_S(\xi)}\right|^{p_*} < \infty. \qquad (9)$$

Let

$$Q_S(y) = \frac{2M_S(y)}{\sigma(y)f_S(y)}.$$

**Theorem 2.** *Suppose that conditions (5) and (9) for some $p_* > 2$ are fulfilled, that the Fisher information $\mathrm{I}(S)$ is continuous at $S(\cdot) = S_*(\cdot)$ and that the law of large numbers holds uniformly in $S(\cdot) \in V_\delta$, that is, for each $\kappa > 0$,*

$$\lim_{T\to\infty}\sup_{S(\cdot)\in V_\delta} \mathbf{P}_S^{(T)}\left\{\left|\frac{1}{T}\int_0^T Q_S(X_t)^2\,dt - \mathbf{E}_S Q_S(\xi)^2\right| > \kappa\right\} = 0. \qquad (10)$$

*Then the empirical moment $\vartheta_T^*$ is an asymptotically efficient estimator of the parameter $\vartheta_S$ under the loss function $\ell(u) = |u|^p$ with $p < p_*$.*

**Proof.** Using the Itô formula, we rewrite the difference $\eta_T = T^{1/2}(\vartheta_T^* - \vartheta_S)$ with a stochastic integral as

$$\eta_T = \frac{H_S(X_T) - H_S(X_0)}{\sqrt{T}} - \frac{1}{\sqrt{T}} \int_0^T Q_S(X_t) \, dW_t \tag{11}$$

under $P_S$. The uniform version of the central limit theorem for stochastic integral (Kutoyants [3], Theorem 3.3.7) allows us to show the weak convergence

$$\mathcal{L}_S \left\{ \frac{1}{\sqrt{T}} \int_0^T Q_S(X_t) \, dW_t \right\} \Longrightarrow \mathcal{N}(0, \mathrm{I}(S)^{-1})$$

uniform in $S(\cdot) \in V_\delta$. Therefore, the empirical estimator is uniformly asymptotically normal

$$\mathcal{L}_S \{ T^{1/2}(\vartheta_T^* - \vartheta_S) \} \Longrightarrow \mathcal{N}(0, \mathrm{I}(S)^{-1}). \tag{12}$$

For the polynomial loss functions we have to verify the uniform integrability of the family of random variables $\{|\eta_T|^p, T \to \infty\}$, but it follows from the representation (11) and the conditions (9) as it was done in Kutoyants [5]. $\square$

**Example 1.** Let us consider the case where $\sigma(y) \equiv 1$ and, for some $\rho > 0, A > 0$ and $L > 0$, define the class of functions

$$\mathbf{S}(\varrho, A, L) = \{ S(\cdot) : |S(y) - S(z)| \leq L|y - z| \text{ for } |y|, |z| \leq A,$$
$$\mathrm{sgn}(y) S(y) \leq -\varrho, \text{ for } |y| \geq A \}.$$

Then, for the function $F(y) = |y|^k$ with any $k > 0$, we can verify all the conditions of Theorem 2 to show that the estimator

$$\vartheta_T^* = \frac{1}{T} \int_0^T |X_t|^k \, dt$$

is consistent, uniformly (in $S(\cdot) \in V_\delta$ with $\delta < \gamma$) asymptotically normal and asymptotically efficient for the polynomial loss functions $\ell(u) = |u|^p$ with any $p > 0$. The verification is quite close to the one given in Kutoyants [5].

**Remark 1.** If we put $F(y) = \chi_{\{y < x\}}$, then $\vartheta = D(x) = \mathbf{P}\{\xi < x\}$ is the value of the invariant distribution function at point $x$. Theorems 1 and 2 yield the asymptotic efficiency of the empirical distribution function

$$\hat{D}_T(x) = \frac{1}{T} \int_0^T \chi_{\{X_t < x\}} \, dt; \tag{13}$$

see Kutoyants [4] for details.

## 3. Parametric estimation

### 3.1. Maximum likelihood estimator

Below we consider the problem of estimation of the drift function of the observed diffusion process described by

$$\mathrm{d}X_t = S(\gamma, X_t)\,\mathrm{d}t + \sigma(X_t)\,\mathrm{d}W_t, \qquad X_0 = x,\ 0 \leq t \leq T. \tag{14}$$

The trend coefficient $S(\cdot)$ is supposed to be known up to an unknown parameter $\gamma \in \Gamma = (\alpha, \beta)$. We suppose that the function $\sigma(\cdot)$ is known and positive, that the conditions (2) and (3) are fulfilled for all $S(\cdot) = S(\gamma, \cdot)$, $\gamma \in \Gamma$, and that the equation (14) has a unique weak solution. Therefore, the process $X_t, t \geq 0$, has an ergodic property and the invariant density function is

$$f(\gamma, y) = G(\gamma)^{-1}\sigma(y)^{-2}\exp\left\{2\int_0^y \frac{S(\gamma, v)}{\sigma(v)^2}\,\mathrm{d}v\right\}$$

with $G(\gamma) = G(S(\gamma, \cdot))$.

The parameter $\vartheta$ is now a function of $\gamma$:

$$\vartheta(\gamma) = \int_{-\infty}^{\infty} F(y)f(\gamma, y)\,\mathrm{d}y, \qquad \gamma \in \Gamma.$$

Set $\Theta = \{\vartheta : \vartheta = \vartheta(\gamma), \gamma \in \Gamma\}$ and denote

$$\dot{\vartheta}(\gamma) = 2\mathbf{E}_\gamma\left[F(\xi)\int_\zeta^\xi \frac{\dot{S}(\gamma, v)}{\sigma(v)^2}\,\mathrm{d}v\right],$$

the derivative of $\vartheta(\gamma)$ with respect to $\gamma$. Here $\xi$ and $\zeta$ are two independent random variables with the same density function $f(\gamma, \cdot)$ and the dot means derivative with respect to $\gamma$.

The regularity condition will be the following:

*Condition C1.* The function $S(\gamma, x)$, $x \in R$, $\gamma \in \Gamma$, has two continuous bounded derivatives on $\gamma$, the function $\sigma(v)^2 \geq \kappa_1 > 0$ for some $\kappa_1$, the Fisher information

$$\mathrm{I}(\gamma) = \int_{-\infty}^{\infty} \left(\frac{\dot{S}(\gamma, y)}{\sigma(y)}\right)^2 f(\gamma, y)\,\mathrm{d}y$$

*is uniformly positive,* $\inf_{\gamma \in \Gamma} \mathrm{I}(\gamma) > 0$, *and the derivative of the function* $\vartheta(\gamma)$ *is separated from zero,*

$$0 < \inf_{\gamma \in \Gamma} |\dot{\vartheta}(\gamma)| \leq \sup_{\gamma \in \Gamma} |\dot{\vartheta}(\gamma)| < \infty. \tag{15}$$

Under this regularity condition, the family of measures $\{\mathbf{P}_\gamma^{(T)}, \gamma \in \Gamma\}$ is LAN at any point $\gamma_0 \in \Gamma$ and we have the Hajek–LeCam lower bound on the risks of all estimators for the loss functions $\ell(u) = |u|^p, p \geq 1$, that is,

$$\lim_{\delta \to 0} \varliminf_{T \to \infty} \inf_{\bar{\vartheta}_T} \sup_{|\gamma - \gamma_0| < \delta} \mathbf{E}_\gamma[\ell(T^{1/2}(\bar{\gamma}_T - \gamma))] \geq \mathbf{E}[\ell(\eta \mathrm{I}(\gamma_0)^{-1/2})], \qquad (16)$$

where $\mathcal{L}(\eta) = \mathcal{N}(0,1)$ (see Kutoyants [3], Theorems 3.3.8 and 3.3.4).

Therefore, the asymptotically efficient estimator in this parametric estimation problem will be defined as follows:

**Definition 2.** *We say that an estimator $\bar{\gamma}_T$ is asymptotically efficient for the loss function $\ell(\cdot)$ if the equality*

$$\lim_{\delta \to 0} \lim_{T \to \infty} \sup_{|\gamma - \gamma_0| < \delta} \mathbf{E}_\gamma \ell(T^{1/2}(\bar{\gamma}_T - \gamma)) = \mathbf{E}\ell(\eta \mathrm{I}(\gamma_0)^{-1/2}) \qquad (17)$$

*holds for all $\gamma_0 \in \Gamma$.*

By condition (15), the above optimality is equivalent to the optimality of an estimator for $\vartheta(\gamma)$: let us put $\mathbf{I}_{\vartheta_0} = \dot{\vartheta}(\gamma_0)^{-2} \mathrm{I}(\gamma_0)$.

**Definition 3.** *The estimator $\bar{\vartheta}_T$ is asymptotically efficient for the loss function $\ell(\cdot)$ if*

$$\lim_{\delta \to 0} \lim_{T \to \infty} \sup_{|\gamma - \gamma_0| < \delta} \mathbf{E}_S \ell(T^{1/2}(\bar{\vartheta}_T - \vartheta(\gamma))) = \mathbf{E}\ell(\eta \mathbf{I}_{\vartheta_0}^{-1/2}). \qquad (18)$$

The condition $\dot{\vartheta}(\gamma) \neq 0$ yields one-to-one mapping $\Theta \leftrightarrow \Gamma$ and the maximum likelihood estimator of the parameter $\vartheta$ is $\hat{\vartheta}_T = \vartheta(\hat{\gamma}_T)$, where $\hat{\gamma}_T$ is the MLE of the parameter $\gamma$ defined by the equation

$$L(\hat{\gamma}_T, \gamma_1, X^T) = \sup_{\gamma \in \Gamma} L(\gamma, \gamma_1, X^T). \qquad (19)$$

Here $\gamma_1$ is some fixed value and the conditional likelihood ratio is given by the formula (see Liptser and Shiryayev [7])

$$\begin{aligned} L(\gamma, \gamma_1, X^T) &\equiv \frac{\mathrm{d}\mathbf{P}_\gamma^{(T)}}{\mathrm{d}\mathbf{P}_{\gamma_1}^{(T)}}(X^T) \\ &= \exp\biggl\{ \int_0^T \frac{S(\gamma, X_t) - S(\gamma_1, X_t)}{\sigma(X_t)^2} \,\mathrm{d}X_t \\ &\quad - \frac{1}{2} \int_0^T [S(\gamma, X_t)^2 - S(\gamma_1, X_t)^2] \sigma(X_t)^{-2} \,\mathrm{d}t \biggr\}. \end{aligned} \qquad (20)$$

It is known that the MLE $\hat{\gamma}_T$ is uniformly consistent, asymptotically normal

$$\mathcal{L}_\gamma\{T^{1/2}(\hat{\gamma}_T - \gamma)\} \Longrightarrow \mathcal{N}(0, \mathrm{I}(\gamma)^{-1})$$

and asymptotically efficient for the polynomial loss function (see Kutoyants [5]).

These properties of $\hat{\gamma}_T$ immediately give the consistency, asymptotic normality

$$\mathcal{L}_\gamma\{T^{1/2}(\hat{\vartheta}_T - \vartheta(\gamma))\} \Longrightarrow \mathcal{N}(0, \mathbf{I}_\vartheta^{-1})$$

and asymptotic efficiency of the MLE $\hat{\vartheta}_T = \vartheta(\hat{\gamma}_T)$.

Therefore, this approach gives us the asymptotically efficient estimator of the moment $\vartheta$, but it has the following disadvantages. The calculation of the MLE $\hat{\gamma}_T$ according to its definition (19) and (20) requires the calculation of stochastic integrals, which are in general not continuous with respect to the uniform topology, as well as maximization of certain nonlinear functional of observations.

## 3.2. One-step estimator based on the empirical estimator

In this subsection, we propose another estimator, which is much easier to calculate and nevertheless is asymptotically efficient in the sense of Definition 2 for the polynomial loss functions.

Note that the empirical estimator

$$\vartheta_T^* = \frac{1}{T}\int_0^T F(X_t)\,\mathrm{d}t$$

given in Section 1 is consistent (with probability 1) and is asymptotically normal

$$\mathcal{L}_\gamma\{T^{1/2}(\vartheta_T^* - \vartheta(\gamma))\} \Longrightarrow \mathcal{N}(0, \mathrm{I}(S(\gamma, \cdot))^{-1}),$$

where $\mathrm{I}(S)$ is defined in (6). Below we will improve this estimator to an asymptotically efficient one. We suppose as well that the equation $\vartheta(\gamma) = \vartheta$ can be solved with respect to $\gamma$ and we have the function $\gamma = \gamma(\vartheta)$ too. This equation can be solved preliminarily, say, numerically because it does not depend on observations and by the condition (15) the solution always exists.

For each locally integrable function $a : \mathbf{R} \to \mathbf{R}$, define

$$G_a(x, \gamma) = -\int_0^x G(S(\gamma, \cdot)) p_{S(\gamma, \cdot)}(y)\left(\int_y^\infty 2a(v) f_{S(\gamma, \cdot)}(v)\,\mathrm{d}v\right)\mathrm{d}y,$$

where

$$p_{S(\gamma, \cdot)}(y) = \exp\left(-2\int_0^y \sigma(y)^{-2} S(\gamma, v)\,\mathrm{d}v\right)$$

if $\int_0^\infty |a(v)| f_{S(\gamma,\cdot)}(v) \, dv < \infty$ (cf. Kutoyants [5] or Yoshida [13] for the estimate of the Green function). Moreover, define a family of functions

$$\mathcal{C}_\gamma = \left\{ a \in C_\uparrow(\mathbf{R}) : \int_\mathbf{R} a(x) f_{S(\gamma,\cdot)}(x) \, dx = 0, [a], G_a(\cdot,\gamma) \in C_\uparrow(\mathbf{R}) \right\},$$

where

$$[a] = -\sigma \nabla G_{a - \langle a, f_{S(\gamma,\cdot)} \rangle}(\cdot, \gamma), \qquad \langle a, f_{S(\gamma,\cdot)} \rangle = \int a(x) f_{S(\gamma,\cdot)}(x) \, dx$$

and $C_\uparrow(\mathbf{R})$ is the space of continuous functions of at most polynomial growth.

Put

$$\Delta_T(\gamma, X^T) = \frac{1}{\sqrt{T}} \int_0^T (\dot{S}(\gamma, X_t) \sigma'(X_t) \sigma(X_t) \\
- \dot{S}(\gamma, X_t) S(\gamma, X_t) - \tfrac{1}{2} \dot{S}'(\gamma, X_t) \sigma(X_t)^2) \sigma(X_t)^{-2} \, dt \qquad (21)$$

and define the estimator

$$\tilde{\gamma}_T = \gamma_T^* + \frac{\Delta_T(\gamma_T^*, X^T)}{\mathrm{I}(\gamma_T^*) \sqrt{T}},$$

where $\gamma_T^* = \gamma(\vartheta_T^*)$. This is the so-called one-step MLE introduced for LAN families by LeCam [6]. It improves a consistent estimator to an efficient one. Note that the one-step estimator

$$\tilde{\vartheta}_T = \vartheta_T^* + \frac{\dot{\vartheta}(\gamma_T^*) \Delta_T(\gamma_T^*, X^T)}{\mathrm{I}(\gamma_T^*) \sqrt{T}}$$

coincides with $\vartheta(\tilde{\gamma}_T)$ up to first order. Thus, first we will consider the former estimator in the sequel. However, $\tilde{\gamma}_T$ is constructed based on our empirical estimator and it is different from the one-step estimator treated in Kutoyants [5].

We will use the following assumptions:

*Condition C2.* (i) *There exists a constant $C$ such that for $\mathcal{H}(x,\gamma) = \partial_\gamma^j \partial_x^i S(\gamma, x)$, $0 \leq i \leq 1, 0 \leq j \leq 3$ and $\sigma(x)$,*

$$\sup_\gamma |\mathcal{H}(x,\gamma)| \leq C(1 + |x|^C)$$

*for all $x \in \mathbf{R}$.*

(ii) *$S(\gamma, \cdot) \in \mathbf{S}(\rho, A, L)$ for some constants $\rho, A,$ and $L$.*

(iii) *For $a_\gamma := \dot{S}(\gamma, \cdot)^2 \sigma(\cdot)^{-2} - \mathrm{I}(\gamma)$, $a_\gamma \in \mathcal{C}_\gamma$ for every $\gamma$ and*

$$\sup_\gamma \mathbf{E}_\gamma[|G_{a_\gamma}(\xi)|] + \sup_\gamma \mathbf{E}_\gamma[[a_\gamma]^2(\xi)] < \infty.$$

**Theorem 3.** *Suppose that* C1 *and* C2 *are fulfilled. Then the estimator* $\tilde{\gamma}_T$ *based on the initial empirical estimator is uniformly consistent, asymptotically normal*

$$\mathcal{L}_\gamma\{T^{1/2}(\tilde{\gamma}_T - \gamma)\} \Longrightarrow \mathcal{N}(0, \mathrm{I}(\gamma)^{-1})$$

*and asymptotically efficient for the polynomial loss function. In particular, one-step estimator $\tilde{\vartheta}_T$ based on the empirical initial estimator has the same asymptotic properties with asymptotic normality:*

$$\mathcal{L}_\gamma\{\sqrt{T}(\tilde{\vartheta}_T - \vartheta(\gamma))\} \Longrightarrow \mathcal{N}(0, \mathbf{I}_\vartheta^{-1}).$$

*Remark 2.* The result may seem to be a corollary to a known general result. However, $L^p$ integrability of the estimator comes from a particular property of the initial estimator such as the empirical estimator we adopted here. Moreover, even if they have the same first-order asymptotic property, one-step estimators are different from each other, depending on the choice of the initial estimator. It will be understood more clearly if we consider the second-order asymptotics.

**Proof of Theorem 3.** Denote the true value of $\gamma$ by $\gamma_0$ and use $\gamma$ as a variable. Define a random field

$$\bar{\Delta}_T(\gamma, X^T) = \frac{1}{\sqrt{T}} \int_0^T \frac{\dot{S}(\gamma, X_t)}{\sigma(X_t)^2} [\mathrm{d}X_t - S(\gamma, X_t)\,\mathrm{d}t].$$

Then, under $P_{\gamma_0}^{(T)}$,

$$\begin{aligned}
\bar{\Delta}_T(\gamma, X^T) &= \frac{1}{\sqrt{T}} \int_0^T \frac{\dot{S}(\gamma, X_t)}{\sigma(X_t)} \mathrm{d}W_t \\
&\quad - \frac{1}{\sqrt{T}} \int_0^T \frac{\dot{S}(\gamma, X_t)}{\sigma(X_t)^2} [S(\gamma, X_t) - S(\gamma_0, X_t)]\,\mathrm{d}t \\
&= \frac{1}{\sqrt{T}} \int_{X_0}^{X_T} \frac{\dot{S}(\gamma, y)}{\sigma(y)^2} \mathrm{d}y - \frac{1}{\sqrt{T}} \int_0^T \frac{\dot{S}(\gamma, X_t)}{\sigma(X_t)^2} S(\gamma_0, X_t)\,\mathrm{d}t \\
&\quad - \frac{1}{2\sqrt{T}} \int_0^T \partial_x \left[\frac{\dot{S}(\gamma, x)}{\sigma(x)^2}\right]\bigg|_{x=X_t} \sigma(X_t)^2\,\mathrm{d}t \\
&\quad - \frac{1}{\sqrt{T}} \int_0^T \frac{\dot{S}(\gamma, X_t)}{\sigma(X_t)^2} [S(\gamma, X_t) - S(\gamma_0, X_t)]\,\mathrm{d}t \\
&= \frac{1}{\sqrt{T}} \int_{X_0}^{X_T} \frac{\dot{S}(\gamma, y)}{\sigma(y)^2}\,\mathrm{d}y \\
&\quad - \frac{1}{\sqrt{T}} \int_0^T \left[s(\gamma, X_t) + \frac{\dot{S}(\gamma, X_t) S(\gamma, X_t)}{\sigma(X_t)^2}\right]\mathrm{d}t,
\end{aligned}$$

where

$$s(\gamma, x) = \frac{1}{2}\sigma(x)^2 \partial_x \left[\frac{\dot{S}(\gamma, x)}{\sigma(x)^2}\right].$$

Given that

$$s(\gamma, x) + \frac{\dot{S}(\gamma, x) S(\gamma, x)}{\sigma(x)^2} = -\frac{1}{\sigma(x)^2}\bigg\{\dot{S}(\gamma, x)\sigma'(x)\sigma(x)$$
$$-\frac{1}{2}\dot{S}'(\gamma, x)\sigma(x)^2 - \dot{S}(\gamma, x)S(\gamma, x)\bigg\},$$

we have

$$\bar{\Delta}_T(\gamma, X^T) = \frac{1}{\sqrt{T}}\int_{X_0}^{X_T}\frac{\dot{S}(\gamma, y)}{\sigma(y)^2}\,dy + \Delta_T(\gamma, X^T). \tag{22}$$

The right-hand side of (22) does not involve Itô stochastic integrals, so it provides a *smooth version* of the random field $\bar{\Delta}_T(\gamma, X^T)$.

Put

$$p(x, \gamma) = \int_0^x \frac{\dot{S}(\gamma, y)}{\sigma(y)^2}\,dy.$$

It is easily seen that under $P_{\gamma_0}^{(T)}$,

$$\sqrt{T}(\tilde{\gamma}_T - \gamma_0) = \sqrt{T}(\gamma_T^* - \gamma_0)$$
$$+ \frac{\Delta_T(\gamma_0, X^T)}{I(\gamma^*)} + \frac{1}{I(\gamma_T^*)}[\Delta_T(\gamma_T^*, X^T) - \Delta_T(\gamma_0, X^T)]$$
$$= \frac{1}{I(\gamma_T^*)\sqrt{T}}\int_0^T \frac{\dot{S}(\gamma_0, X_t)}{\sigma(X_t)}\,dW_t + R_T(\gamma_0),$$

where $R_T(\gamma_0) = R_{1,T}(\gamma_0) + R_{2,T}(\gamma_0)$ with

$$R_{1,T}(\gamma_0) = -\frac{1}{I(\gamma_T^*)\sqrt{T}}[p(X_T, \gamma_0) - p(X_0, \gamma_0)]$$

and

$$R_{2,T}(\gamma_0) = \sqrt{T}(\gamma_T^* - \gamma_0) + \frac{1}{I(\gamma_T^*)}[\Delta_T(\gamma_T^*, X^T) - \Delta_T(\gamma_0, X^T)].$$

It follows from Condition C2 that for $T > 1$,

$$\sup_{\gamma_0} \mathbf{E}_{\gamma_0}[|R_{1,T}(\gamma_0)|] \leq \frac{C}{\sqrt{T}}.$$

Next, we will estimate the second residual $R_{2,T}(\gamma_0)$. Denote by $L_\gamma$ the generator of the diffusion process $X_t$ that corresponds to the parameter $\gamma$,

$$L_\gamma = S(\gamma, x)\frac{\mathrm{d}}{\mathrm{d}x} + \frac{1}{2}\sigma(x)^2\frac{\mathrm{d}^2}{\mathrm{d}x^2},$$

and put $M(\gamma, x) = \partial_\gamma(L_\gamma p)(x, \gamma)$. Then obviously $M(\gamma, x) = L_\gamma \dot{p}(x, \gamma) + \dot{S}(\gamma, x)^2 \sigma(x)^{-2}$. We have

$$R_{2,T}(\gamma_0) = \frac{R_{3,T}(\gamma_0)}{\mathrm{I}(\gamma_T^*)}\sqrt{T}(\gamma_T^* - \gamma_0),$$

where

$$R_{3,T}(\gamma_0) = -\frac{1}{T}\int_0^T M(\gamma_T^{**}, X_t)\,\mathrm{d}t + \mathrm{I}(\gamma_T^*).$$

Then $R_{3,T}(\gamma_0) = \sum_{i=4}^8 R_{i,T}(\gamma_0)$ with

$$R_{4,T}(\gamma_0) = \mathrm{I}(\gamma_T^*) - \mathrm{I}(\gamma_0),$$

$$R_{5,T}(\gamma_0) = -\frac{1}{T}\int_0^T \dot{S}(\gamma_0, X_t)^2 \sigma(X_t)^{-2}\,\mathrm{d}t + \mathrm{I}(\gamma_0),$$

$$R_{6,T}(\gamma_0) = -\frac{1}{T}\int_0^T [\dot{S}(\gamma_T^{**}, X_t)^2 - \dot{S}(\gamma_0, X_t)^2]\sigma(X_t)^{-2}\,\mathrm{d}t,$$

$$R_{7,T}(\gamma_0) = -\frac{1}{T}\int_0^T (L_{\gamma_0}\dot{p})(\gamma_0, X_t)\,\mathrm{d}t,$$

and

$$R_{8,T}(\gamma_0) = -\frac{1}{T}\int_0^T [(L_{\gamma_T^{**}}\dot{p})(\gamma_T^{**}, X_t) - (L_{\gamma_0}\dot{p})(\gamma_0, X_t)]\,\mathrm{d}t.$$

By using the uniform non-degeneracy of $\dot{\vartheta}(\gamma)$ (Condition C1) and the uniform estimate for $\vartheta_T^*$, we obtain, for every $p > 1$,

$$\sup_{\gamma_0} \mathbf{E}_{\gamma_0}[|\sqrt{T}(\gamma_T^* - \gamma_0)|^p] \leq \left(\sup_\gamma |\dot{\vartheta}(\gamma)|\right)^p \cdot \sup_{\gamma_0} \mathbf{E}_{\gamma_0}[|\sqrt{T}(\vartheta_T^* - \vartheta_0)|^p] \leq C_p < \infty$$

for all $T \in \mathbf{R}_+$. Because $\sup_\gamma |\dot{\mathrm{I}}(\gamma)| < \infty$ as a consequence of Condition C2, we obtain for $i = 4$,

$$\sup_{\gamma_0} \mathbf{E}_{\gamma_0}[|R_{i,T}(\gamma_0)|] \leq \frac{C}{\sqrt{T}}. \tag{23}$$

In a similar fashion, we obtain the estimate (23) for $i = 6, 8$. Also, $R_{5,T}$ and $R_{7,T}$ can be estimated by using Itô's formula, the Burkholder–Davis–Gundy inequality and Condition

C2. Thus, $R_T(\gamma_0)$ is estimated uniformly in $\gamma_0$. Finally, the uniform central limit theorem for the principal term of the $\sqrt{T}(\tilde{\gamma}_T - \gamma_0)$ implies the desired result.

The assertion for $\tilde{\vartheta}$ is an easy consequence. □

**Remark 3.** If we are interested in estimating the distribution function [see Remark 1 with $\vartheta = D(\gamma, x)$] and if $\hat{D}_T(x)$ is its empirical distribution function (13), then according to Theorem 3, the estimator

$$\tilde{D}_T(x) = \hat{D}_T(x) + \frac{\dot{D}(\gamma_T^*, x)\Delta_T(\gamma_T^*, X^T)}{\mathrm{I}(\gamma_T^*)\sqrt{T}}$$

will be asymptotically efficient in the following sense:

$$\lim_{\delta \to 0} \varliminf_{T \to \infty} \sup_{|\gamma - \gamma_0| < \delta} \mathbf{E}_\gamma \ell(T^{1/2}(\tilde{D}_T(x) - D(\gamma, x))) = \mathbf{E}\ell(\eta \mathbf{I}_{D_0}^{-1/2}).$$

Here $\gamma_T^*$ is the solution of the equation $D(\gamma_T^*, x) = \hat{D}_T(x)$, and $\mathrm{I}(\gamma_0)$ and $\mathbf{I}_{D_0}$ are the corresponding Fisher informations (see Kutoyants [5]).

## 4. Asymptotic expansion

We are still considering the diffusion process $X_t$ that satisfies the stochastic differential equation (1) and the nonparametric estimator for the expectation $\vartheta = \mathbf{E}[F(\xi)]$. As shown in the previous sections,

$$\vartheta_T^* = \frac{1}{T}\int_0^T F(X_t)\,\mathrm{d}t$$

is asymptotically normal and asymptotically efficient in a nonparametric sense. In this section, we shall investigate the higher order distribution of our nonparametric estimator, more precisely, the asymptotic expansion of its distribution will be presented.

The asymptotic expansion of the distribution of ergodic diffusions recently was obtained by Yoshida [12, 13, 14] applying the Malliavin calculus. Among two possible methodologies, that is, the *global approach* and the *local approach*, here we shall take the newly developed local approach formalized by Kusuoka and Yoshida [2] for continuous-time processes and applied by Sakamoto and Yoshida [8]. The support theorems serve to verify the non-degeneracy (see Yoshida [14]).

Let $C_\uparrow^\infty(\mathbf{R})$ be the space of smooth functions on $\mathbf{R}$, all derivatives of which are of at most polynomial growth, and let $C_B^\infty(\mathbf{R})$ be the space of bounded smooth functions on $\mathbf{R}$ with bounded derivatives. We denote by $B\mathcal{G}$ the set of bounded $\mathcal{G}$-measurable functions. We assume that $F, S \in C_\uparrow^\infty(\mathbf{R})$, $S', \sigma \in C_B^\infty(\mathbf{R})$, and $\sigma(x) > 0$ for any $x \in \mathbf{R}$, and that $X_t$ is stationary. We may construct a solution $X_t$ over a (partial) Wiener space $(\Omega, P)$, that is, $\Omega = \mathbf{R} \times W$, where $W = \{w\colon \mathbf{R}_+ \to \mathbf{R}, \text{continuous } w(0) = 0\}$, and $P = \nu \otimes \tilde{P}$, $\tilde{P}$ being a Wiener measure on $W$ and $\nu$ being the stationary distribution of the diffusion process. Let

$$\mathcal{B}_I^X = \sigma[X_t : t \in I] \vee \mathcal{N}$$

for $I \in \mathbf{R}_+$, $\mathcal{N}$ being the $\sigma$-field generated by null sets. We later use the following conditions:

Condition A1. There exists a positive constant $a$ such that

$$\|\mathbf{E}[h|\mathcal{B}_{[s]}^X] - \mathbf{E}[h]\|_1 \leq a^{-1} e^{-a(t-s)} \|h\|_\infty$$

for any $s, t \in \mathbf{R}_+$, $s \leq t$, and for any $h \in B\mathcal{B}_{[t,\infty)}^X$.

Condition A2. Process $X_0 \in \bigcap_{p>1} L^p(P)$.

A sufficient condition for A1 is provided, for example, in Veretennikov [10, 11] and Kusuoka and Yoshida [2]. Indeed, if $\sigma \in C_B^\infty(\mathbf{R})$ and if there exists a function $\rho \in C^\infty(\mathbf{R})$ such that $\rho > 0$, $\int_\mathbf{R} \rho(x)\,\mathrm{d}x = 1$ and $\limsup_{|x|\to\infty} \rho(x)^{-1} L^* \rho(x) < 0$, where $L^*$ is the formal adjoint operator of the generator $L$ of this diffusion process, then Condition A1 holds true (see Kusuoka and Yoshida [2]).

Put

$$Z_T = \int_0^T q(X_t)\,\mathrm{d}t,$$

where $q(x) = F(x) - \vartheta$. As in Kusuoka and Yoshida [2], we define the $r$th cumulant function of $Z_T/\sqrt{T}$ by

$$\chi_{T,r}(u) = \left(\frac{\mathrm{d}}{\mathrm{d}\varepsilon}\right)_0^r \log \mathbf{E}\left[\exp\left(i\varepsilon u \cdot \frac{Z_T}{\sqrt{T}}\right)\right].$$

Then define function

$$\hat{\Psi}_{T,k}(u) = \exp(\tfrac{1}{2}\chi_{T,2}(u)) + \sum_{r=1}^k T^{-r/2} \tilde{P}_{T,r}(u),$$

where $\tilde{P}_{T,r}(u)$ are defined by the formal Taylor expansion

$$\exp\left(\sum_{r=2}^\infty \frac{1}{r!}\varepsilon^{r-2}\chi_{T,r}(u)\right) = \exp\left(\frac{1}{2}\chi_{T,2}(u)\right) + \sum_{r=1}^\infty \varepsilon^r T^{-r/2} \tilde{P}_{T,r}(u).$$

Denote by $\Psi_{T,k}$ the signed measure defined as the Fourier inversion of $\hat{\Psi}_{T,k}$. For measurable function $h: \mathbf{R} \to \mathbf{R}$, let

$$\omega(h,r) = \int_\mathbf{R} \sup\{|h(x+y) - h(x)| : |y| \leq r\} \phi(x; 0, \mathrm{I}_\star^{-1})\,\mathrm{d}x,$$

where $\phi(x; \mu, \Sigma)$ denotes the probability density of the normal distribution $N(\mu, \Sigma)$ and $\mathrm{I}_\star$ is an arbitrary positive number smaller than $I_*$. The Hermite polynomials are defined

by
$$h_k(z;\Sigma) = (-1)^k \phi(z;0,\Sigma)^{-1} \partial_z^k \phi(z;0,\Sigma)$$

for a positive constant $\Sigma$.

For a continuous function $a: \mathbf{R} \to \mathbf{R}$, we define $G_a : \mathbf{R} \to \mathbf{R}$ by
$$G_a(x) = \int_{-\infty}^{x} G(S) p(y) \left( \int_{-\infty}^{y} 2a(v) f(v) \, \mathrm{d}v \right) \mathrm{d}y,$$

where $p(y) = \exp(-2 \int_0^y \sigma(v)^{-2} S(v) \, \mathrm{d}v)$ if the mapping
$$y \mapsto p(y) \left( \int_{-\infty}^{y} a(v) f(v) \, \mathrm{d}v \right)$$

is in $L^1((-\infty, 0], \mathrm{d}y)$. We write $\langle a, f \rangle = \int_{\mathbf{R}} a(x) f(x) \, \mathrm{d}x$ and
$$[a] = -\sigma \nabla G_{a - \langle a, f \rangle}.$$

Define the set of functions
$$\mathcal{C} = \left\{ a \in C_\uparrow(\mathbf{R}) \mid \int_{-\infty}^{\infty} a(x) f(x) \, \mathrm{d}x = 0; \right.$$
$$\left. p(\cdot) \int_{-\infty}^{\cdot} a(x) f(x) \, \mathrm{d}x \in L^1((-\infty, 0]); [a], G_a \in C_\uparrow(\mathbf{R}) \right\}.$$

**Theorem 4.** *Let $k \in \mathbf{N}$ and let $M, \gamma, K > 0$. Suppose that $F: \mathbf{R} \to \mathbf{R}$ is not constant. Then*

(1) *There exist constants $\delta > 0$ and $c > 0$ such that for $h \in \mathcal{E}(M, \gamma)$,*
$$|\mathbf{E}[h(\sqrt{T}(\vartheta_T^* - \vartheta))] - \Psi_{T,k}[h]| \leq c\omega(h, T^{-K}) + \varepsilon_T^{(k)},$$

*where $\varepsilon_T^{(k)} = o(T^{-(k+\delta)/2})$. Here $\mathcal{E}(M, \gamma) = \{h: \mathbf{R} \to \mathbf{R}, \text{ measurable, } |h(x)| \leq M(1 + |x|)^\gamma (x \in \mathbf{R})\}$.*

(2) *The signed measure $\mathrm{d}\Psi_{T,1}$ has a density $\mathrm{d}\Psi_{T,1}(z)/\mathrm{d}z = p_{T,1}(z)$ with*
$$p_{T,1}(z) = \phi(z; 0, \kappa_T^{(2)})(1 + \tfrac{1}{6}\kappa_T^{(3)} h_3(z; \kappa_T^{(2)})),$$

*where $\kappa_T^{(r)}$ is the rth cumulant of $\sqrt{T}(\vartheta_T^* - \vartheta)$. Moreover, if $q$ and $[q]^2 - \langle [q]^2, f \rangle$ belong to $\mathcal{C}$, then*
$$p_{T,1}(z) = p_{T,1}^*(z) + R_T(z),$$

*where*
$$p_{T,1}^*(z) = \phi(z; 0, \mathrm{I}_*^{-1}) \left( 1 + \frac{1}{2\sqrt{T}} \mathbf{E}[[[q]^2][q](\xi)] h_3(z; \mathrm{I}_*^{-1}) \right)$$

*and*

$$\lim_{T \to \infty} \sqrt{T} \sup_{z \in \mathbf{R}} \{|R_T(z)| \exp(bz^2)\} = 0$$

*for some positive constant b. In particular,*

$$\left| \mathbf{E}[h(\sqrt{T}(\vartheta_T^* - \vartheta))] - \int_{\mathbf{R}} h(z) p_{T,1}^*(z) \, dz \right| \leq c\omega(h, T^{-K}) + \tilde{\varepsilon}_T$$

*for any $h \in \mathcal{E}(M, \gamma)$, with $\tilde{\varepsilon}_T = o(1/\sqrt{T})$.*

**Proof.** We consider the stochastic flow $\bar{X}(t,x) = (X(t,x), Z(t,x))$ that satisfies the stochastic differential equation

$$d\bar{X}(t,x) = \bar{V}_0(\bar{X}(t,x)) \, dt + \bar{V}(\bar{X}(t,x)) \circ dW_t, \qquad \bar{X}(0,x) = (x,0),$$

where $\bar{V}_0(x,z) = (S(x) - 2^{-1}\sigma'(x)\sigma(x), q(x))$ and $\bar{V}(x,z) = (\sigma(x), 0)$. By assumption, $F$ is not constant; hence, there exists a point $x_0 \in \mathbf{R}$ such that $F'(x) \neq 0$ in a neighborhood $U$ of $x_0$. Easy calculus shows that

$$\mathrm{Span}\{\bar{V}(x_0, 0), [\bar{V}, \bar{V}_0](x_0, 0)\} = \mathbf{R}^2.$$

In the same argument used to prove Theorem 4 of Kusuoka and Yoshida [2] or that in Example 2 of the same paper, Condition [A3'] of it can be verified. Indeed, take the sequence $u(j) = j$, $v(j) = j+1$, $j \in \mathbf{R}$, and let $\psi_j = \varphi(X_j)$ for some truncation function $\varphi \in C^{\infty}(\mathbf{R}; [0,1])$ with compact support in $U$ taking value 1 near $x_0$. For each $j \in \mathbf{N}$, the Malliavin operator $L_j$ is constructed in the usual way so that $L_j$ does not shift the path $w$ outside of $[j, j+1]$. Then it is known that for sufficiently small $U$, $\{(\det \sigma_{\bar{X}(1,x)})^{-1}; x \in U\}$ is bounded in $L^p(P)$ for any $p > 1$. Therefore, Condition [A3'] of Kusuoka and Yoshida [2] can be verified. Thus, the first assertion has been obtained.

For the formula of $p_{T,k}$ with the validity, see Sakamoto and Yoshida [8, 9]. To obtain the last result, some calculations are involved. Given that $q \in \mathcal{C}$,

$$\frac{1}{\sqrt{T}} \int_0^T q(X_t) \, dt = \frac{G_q(X_T) - G_q(X_0)}{\sqrt{T}} + \frac{1}{\sqrt{T}} \int_0^T [q](X_t) \, dW_t \qquad (24)$$

and $\mathbf{E}[|G_q(\xi)|^p] < \infty$ for any $p > 1$. Let $0 < \beta < 1/2$. The strong mixing property A1 induces the so-called covariance inequality, which yields

$$\mathbf{E}\left[G_q(X_T) \cdot \frac{1}{\sqrt{T}} \int_{T^\beta}^{T-T^\beta} [q](X_t) \, dW_t\right]$$

$$\leq 8\alpha(\mathcal{B}_{[0, T-T^\beta]}, \mathcal{B}_{[T, \infty)})^{1/r} \left\| \frac{1}{\sqrt{T}} \int_{T^\beta}^{T-T^\beta} [q](X_t) \, dW_t \right\|_p \|G_q(X_T)\|_q$$

$$\leq a^{-1} \exp\left(-\frac{aT^\beta}{r}\right),$$

where $\alpha$ is the coefficient of the strong mixing, $a$ is some positive constant and $1/p + 1/q + 1/r = 1$. A similar estimate holds even if $G_q(X_T)$ above is replaced with $G_q(X_0)$. Accordingly, we obtain

$$\kappa_T^{(2)} = \frac{1}{T}\int_0^T \mathbf{E}[[q]^2(X_t)]\,dt + \frac{2}{\sqrt{T}}\mathbf{E}\left[(G_q(X_T) - G_q(X_0)) \cdot \frac{1}{\sqrt{T}}\int_0^T [q](X_t)\,dW_t\right]$$

$$+ \frac{1}{T}\mathbf{E}[(G_q(X_T) - G_q(X_0))^2]$$

$$= \mathbf{E}[[q]^2(\xi)] + o\left(\frac{1}{\sqrt{T}}\right). \tag{25}$$

Next, we consider the third cumulant $\kappa_T^{(3)}$. Put $k = [q]$ and denote simply $k_t = k(X_t)$. By assumption, $k \in C_\uparrow(\mathbf{R})$, and so we see by Itô's lemma that

$$\mathbf{E}\left[\left(\frac{1}{\sqrt{T}}\int_0^T k_t\,dW_t\right)^3\right]$$

$$= 3T^{-3/2}\mathbf{E}\left[\int_0^T \left(\int_0^t k_s\,dW_s\right)k_t^2\,dt\right]$$

$$= 3T^{-3/2}\mathbf{E}\left[\int_0^T \left(\int_0^t k_s\,dW_s\right)(k_t^2 - \langle k^2, f\rangle)\,dt\right]$$

$$= 3T^{-1/2}\mathbf{E}\left[\frac{1}{\sqrt{T}}\int_0^T k_s\,dW_s \cdot \frac{1}{\sqrt{T}}\int_0^T (k_t^2 - \langle k^2, f\rangle)\,dt\right].$$

Because $k^2 - \langle k^2, f\rangle \in \mathcal{C}$, the right-hand side equals

$$\frac{3}{\sqrt{T}}\mathbf{E}\left[\frac{1}{\sqrt{T}}\int_0^T k_s\,dW_s\left\{\frac{G_{k^2 - \langle k^2, f\rangle}(X_T) - G_{k^2 - \langle k^2, f\rangle}(X_0)}{\sqrt{T}} + \frac{1}{\sqrt{T}}\int_0^T [[q]^2](X_t)\,dW_t\right\}\right]$$

$$= \frac{3}{\sqrt{T}}\mathbf{E}[[q]^2[q](\xi)] + o\left(\frac{1}{\sqrt{T}}\right). \tag{26}$$

In view of (24) and by a similar argument after it, it is possible to see that the cross terms [e.g., those between $G_q(X_T) - G_q(X_0)$ and $(T^{-1/2}\int_0^T [q]\,dW_t)^2$] have asymptotically no contribution to $\kappa_T^{(3)}$. It follows from (25) and (26) that

$$\sup_{z \in \mathbf{R}} e^{bz^2}|p_{T,1}(z) - p_{T,1}^*(z)| = o\left(\frac{1}{\sqrt{T}}\right)$$

for some positive constant $b$. This completes the proof. $\square$

# Acknowledgement

We would like to thank the referees for helpful comments.